\documentclass[11pt]{article}
\usepackage{epsfig,amssymb}
\usepackage{palatino}
\usepackage{amsmath}

\def \ra {{\quad\Rightarrow\quad}}
\def \lr {{\quad\Leftrightarrow\quad}}

\def\R{{\mathbb R}}
\def\N{{\mathbb N}}
\def\SS{{\cal S}}
\def\FF{{\cal F}}
\def\e{\varepsilon}
\def \l {\ell}
\def \Pn {P_{\Delta_n}}
\def \la {\langle}
\def \r {\rangle}

\def \supp {\mbox{\rm supp\,}}
\def \spa {\mbox{\rm span\,}}
\newtheorem{lemma}{Lemma}[section]
\newtheorem{proposition}[lemma]{Proposition}
\newtheorem{corollary}[lemma]{Corollary}
\newtheorem{theorem}[lemma]{Theorem}
\newtheorem{remark}[lemma]{Remark}

\newtheorem{definition}[lemma]{Definition}

\def\be  {\begin{equation}} 
\def\ee  {\end{equation}}
\def\ba  {\begin{eqnarray}} 
\def\ea  {\end{eqnarray}}
\def\baa {\begin{eqnarray*}} 
\def\eaa {\end{eqnarray*}}
\makeatletter
\@addtoreset{equation}{section}
\makeatother
\def\proof{\medskip\noindent{\bf Proof.} }
\def\qed{\hfill $\Box$}
\newcommand {\lb} {\label}
\newcommand {\rf[1]} {(\ref{#1})}

\begin{document}

\title{On almost everywhere convergence
of orthogonal spline projections with arbitrary knots}

\author{M.\,Passenbrunner, 
A.\,Shadrin }

\maketitle

\abstract{
The main result of this paper is a proof that, for any $f \in L_1[a,b]$, 
a sequence of its orthogonal projections $(P_{\Delta_n}(f))$ onto splines 
of order $k$ with arbitrary knots $\Delta_n$, converges almost everywhere 
provided that the mesh diameter $|\Delta_n|$ 
tends to zero, namely
$$
     f \in L_1[a,b] \ra
     P_{\Delta_n}(f,x) \to f(x) \quad \mbox{a.e.} \quad (|\Delta_n|\to 0)\,.
$$
This extends the earlier result that, 
for $f \in L_p$,  we have convergence
$P_{\Delta_n}(f) \to f$ in the $L_p$-norm for $1 \le p \le \infty$, where we interpret $L_\infty$ as the space of continuous functions.
}


\section{Introduction}


Let an interval $[a,b]$ and $k \in \N$ be fixed. For a knot-sequence 
$\Delta_n = (t_i)_{i=1}^{n+k}$ such that 
\baa
   t_i \le t_{i+1}, && t_i < t_{i+k}, \\
   t_1 = \cdots = t_k = a, && b = t_{n+1} = \cdots = t_{n+k},
\eaa
let $(N_i)_{i=1}^n$ be the sequence of $L_\infty$-normalized B-splines 
of order $k$ on $\Delta_n$ forming a partition of unity, with the properties
$$
   \supp N_i = [t_i, t_{i+k}]\,, \qquad N_i \ge 0\,, \qquad
   \sum_i N_i \equiv 1\,.
$$
For each $\Delta_n$, we define then the space $\SS_k(\Delta_n)$ of splines 
of order $k$ with knots $\Delta_n$ as the linear span of $(N_i)$, namely 
$$
   s \in \SS_k(\Delta_n) \lr s = \sum_{i=1}^n c_i N_i\,,\quad c_i \in \R\,,
$$
so that $\SS_k(\Delta_n)$ is the space of piecewise polynomial functions
of degree $\le k-1$, with $k-1-m_i$ continuous derivatives at $t_i$, 
where $m_i$ is multiplicity of $t_i$. 
Throughout the paper, we use the following notations: 
\baa
   & I_i := [t_i,t_{i+1}], \quad h_i := |I_i| := t_{i+1}-t_i\,, & \\
   & E_i := [t_i,t_{i+k}], \quad \kappa_i := |E_i| := t_{i+k}-t_i\,,&
\eaa
where $E_i$ is the support of the B-spline $N_i$. 
With ${\rm conv}(A,B)$ standing for the convex hull of two sets $A$ and $B$,
we also set   
\baa
&  I_{ij} := {\rm conv}(I_i,I_j)
= [t_{\min(i,j)},t_{\max(i,j) + 1}]\,, & \\
&  E_{ij}  := {\rm conv}(E_i,E_j)
= [t_{\min(i,j)},t_{\max(i,j) + k}]\,. &
\eaa
Finally, $|\Delta_n| := \max_i |I_i|$ is the mesh diameter of $\Delta_n$.

Now, let $P_{\Delta_n}$ be the orthoprojector onto $\SS_k(\Delta_n)$ with 
respect to the ordinary inner product $\la f,g\r = \int_a^b f(x)g(x)\,dx$, 
i.e.,
$$
    \la P_{\Delta_n} f,s\r = \la f,s\r, \qquad 
    \forall s \in \SS_k(\Delta_n)\,,
$$
which is well-defined for $f \in L_1[a,b]$. 

Some time ago, one of us proved \cite{s2} de\,Boor's conjecture 
that the max-norm of $\Pn$ 
is bounded independently of the knot-sequence, i.e.,
\be \lb{dBc}
    \sup_{\Delta_n} \|\Pn\|_\infty < c_k\,.
\ee
This readily implies 
convergence of orthogonal spline projections 
in the $L_p$-norm,
\be \lb{p}
     f \in L_p[a,b] \ra
     P_{\Delta_n}(f) \stackrel{L_p}{\to} f\,, \qquad 
     1 \le p \le \infty\,,
\ee
where we interpret $L_\infty$ as $C$, the space of continuous functions.
In this paper, we prove that the max-norm boundedness 
of $\Pn$ implies also almost everywhere (a.e.) convergence of 
orthogonal projections $(P_{\Delta_n}(f))$ with arbitrary knots $\Delta_n$ 
provided that the mesh diameter $|\Delta_n|$ 
tends to zero. 

The main outcome of this article  is the following statement.

\begin{theorem} \lb{t1}
For any $k \in \N$ and any 
sequence of partitions $(\Delta_n)$ such that $|\Delta_n| \to 0$, we have
\be \lb{Ptof}
    f \in L_1[a,b] \ra P_{\Delta_n}(f,x) \to f(x) \quad \mbox{a.e.}
\ee
\end{theorem}

The proof is based on a standard approach of verifying two conditions
which imply a.e.\ convergence for $f \in L_1$:

1) there is a dense subset $\FF$ of $L_1$ such that $\Pn(f,x) \to f(x)$ a.e.
for $f \in \FF$,

2) the maximal operator $P^*(f,x) := \sup_n |\Pn(f,x)|$ is of the weak 
$(1,1)$-type, 
\be \lb{P}
    m \{x \in [a,b]: P^*(f,x) > t\} < \frac{c_k}{t}\,\|f\|_1,
\ee
with $mA$ being the Lebesgue measure of $A$. 
The first condition is easy: 
by \rf[p], a.e.\ convergence 
(in fact, uniform convergence) takes place for continuous functions,
\be \lb{C}
     f \in C[a,b] \ra
     P_{\Delta_n}(f,x) \to f(x) \quad \mbox{uniformly in $x$.}
\ee
For the non-trivial part \rf[P], we prove a stronger inequality
of independent interest, namely that 
\be \lb{M}
     |\Pn(f,x)| \le c_k M(f,x)\,,
\ee
where $M(f,x)$ is the Hardy--Littlewood maximal function. It
satisfies a weak $(1,1)$-type inequality, hence \rf[P] holds too.

The main technical tool which leads to \rf[M]
is a new estimate for the elements $\{a_{ij}\}$ of the inverse 
of the Gram matrix of the B-spline functions, which reads as follows.

\begin{theorem} \lb{p1}
For any $\Delta_n$, let $\{a_{ij}\}_{i,j=1}^n$ be the inverse of the B-spline 
Gram matrix $\{\langle N_i,N_j\rangle\}$.
Then, 
\be \lb{a0}
    |a_{ij}| \le K \gamma^{|i-j|} h_{ij}^{-1}\,, 
\ee
where 
$$
   h_{ij} := \max \{h_s: I_s \subset E_{ij}\}\,,
$$
and $K>0$ and $\gamma \in (0,1)$ are constants which depend only on $k$, but
not on $\Delta_n$.
\end{theorem}
A pass from \rf[a0] to \rf[M] proceeds as follows.
Let  $K_{\Delta_n}$ be the Dirichlet kernel of the operator $P_{\Delta_n}$,
defined by the relation 
$$
   P_{\Delta_n}(f,x) = \int_a^b K_{\Delta_n}(x,y)f(y)\, dy\,, \quad 
  \forall f \in L_1[a,b].
$$
Then, \rf[a0] implies the inequality
\be \lb{K}
     |K_{\Delta_n}(x,y)| 
\leq C\,\theta^{|i-j|}|I_{ij}|^{-1}, \qquad x\in I_i,\ y\in I_j,
\ee
where $C > 0$ and $\theta \in (0,1)$. 
Now, \rf[M] is immediately obtained from \rf[K].

With a bit more sophisticated arguments, though still standard ones, estimate \rf[K] on $K_{\Delta_n}$ allows us also to prove 
convergence of $\Pn f$ at Lebesgue points of $f$. The latter forms a set 
of full measure, so we derive this refinement of Theorem \ref{t1} as a 
byproduct. 

Estimate \rf[a0] is also useful
in other applications, for instance 
in \cite{p} it is applied to 
obtain unconditionality of orthonormal spline bases with arbitrary 
knot-sequences in $L_p$-spaces for $1 < p < \infty$.

We note that, previously, a.e.\ convergence of spline orthoprojections
was studied by Ciesielski \cite{c1} who established
\rf[Ptof] for dyadic partitions with any $k \in \N$,
and by Ciesielski--Kamont \cite{ck} who proved this result for any $\Delta_n$
with $k=2$, i.e., for linear splines. Both papers used \rf[a0] as an 
intermediate step, however our proof of \rf[a0] for all $k$ with arbitrary 
knots $\Delta_n$ is based on quite different arguments. 
The main difference is that the proof of \rf[a0] for linear 
splines in \cite{ck} does not rely on the 
mesh-independent bound \rf[dBc] for $\|P_{\Delta_n}\|_\infty$, and 
can be used to get such a bound for linear splines, whereas our proof
depends on \rf[dBc] in an essential manner.

The paper is organised as follows. 
In Sect.\,\ref{pt1}, we show how Theorem \ref{p1} leads to \rf[K] and the 
latter to \rf[M]. We complete then the proof of a.e.\ convergence of $(\Pn(f))$ 
using the scheme indicated above.
In Sect.\,\ref{clp}, as a byproduct, we show that $(\Pn(f))$ converges
at Lebesgue 
points, thus characterizing the convergence set in a sense. 
Theorem \ref{p1} is proved then in Sect.\,\ref{pp1}
based on Lemma \ref{le1}, which lists several specific properties
of the inverse $\{a_{ij}\}$ of the B-spline Gram matrix 
$G_0 := \{\langle N_i,N_j\rangle\}$. 
Those properties are proved in 
the final Sect.\,\ref{pl1}, and they are based mostly 
on Demko's theorem on the inverses of band matrices, which we apply
to the rescaled Gram matrix $G := (\la M_i, N_j\r)$, where 
$M_i := \frac{k}{\kappa_i} N_i$. The uniform bound 
$\|G^{-1}\|_\infty < c_k$, being equivalent to \rf[dBc], plays a crucial role here.


\section{Proof of Theorem \ref{t1} } \lb{pt1}


Here, we prove the weak-type inequality \rf[P], 
then  recall a simple proof of \rf[C], and as a result deduce 
the a.e.\ convergence for all $f \in L_1$. 

\medskip
We begin with an estimate for the Dirichlet kernel $K_{\Delta_n}$.
\begin{lemma}
For any $\Delta_n$, the Dirichlet kernel $K_{\Delta_n}$ 
satisfies the inequality
\be \lb{Kineq}
    |K_{\Delta_n}(x,y)| 
\leq C\,\theta^{|i-j|}|I_{ij}|^{-1},\qquad x\in I_i,\,y\in I_j,
\ee
where $C > 0$ and $\theta \in (0,1)$ 
are constants that depends only on $k$.
\end{lemma}
\proof
First note that, with the inverse $\{a_{\ell m}\}$ of the 
B-spline Gram matrix $\{\langle N_\ell,N_m\rangle\}$,
the Dirichlet kernel $K_{\Delta_n}$ can be written in the form
$$
    K_{\Delta_n}(x,y) = \sum_{\ell,m=1}^n a_{\ell m}N_\ell(x)N_m(y).
$$
For $x\in I_i$ and $y\in I_j$, since $\supp N_\ell = [t_\ell,t_{\ell+k}]$ 
and $\sum N_\ell(x)N_m(y) \equiv 1$,  
we obtain
$$
   |K_{\Delta_n}(x,y)| 
\le \max_{\substack{i-k+1\leq \ell\leq i \\ j-k+1\leq m\leq j}} |a_{\ell m}|\,.
$$
Next, we rewrite inequality \rf[a0] for $a_{\ell m}$ 
in terms of  $E_{\ell m} = [t_{\min(\ell,m)}, t_{\max(\ell,m) + k}]$: 
as  $h_{\ell m}$ is the largest knot-interval in $E_{\ell m}$, 
we have $h_{\ell m}^{-1} \leq (|\ell-m|+k) |E_{\ell m}|^{-1}$, hence 
for any real number $\theta\in(\gamma,1)$,
$$
     |a_{\ell m}|
\le K \gamma^{|\ell-m|} (|\ell-m| + k)|E_{\ell m}|^{-1}
\le C_1 \theta^{|\ell-m|}|E_{\ell m}|^{-1},
$$
where $C_1$ depends on $k$ and $\theta$. 
Therefore, 
$$
     |K_{\Delta_n}(x,y)| 
\le C_1 \max_{\substack{i-k+1\leq \ell\leq i \\ j-k+1\leq m\leq j}} 
      \theta^{|\ell-m|}|E_{\ell m}|^{-1}\,.
$$
For indices $\ell$ and $m$ in the above maximum, we have
$I_{ij}\subset E_{\ell m}$, hence $|E_{\ell m}|^{-1} \le |I_{ij}|^{-1}$, and
also $|\l-m| > |i-j|-k$, hence $\theta^{|\ell-m|} \le \theta^{-k} \theta^{|i-j|}$, 
and inequality \rf[Kineq] follows.
\qed

\begin{definition} \rm
For an integrable $f$, the Hardy--Littlewood maximal 
function is defined as
\be \lb{M0}
    M(f,x) := \sup_{I \ni x} |I|^{-1} \int_I |f(t)|\,dt\,,
\ee
with the supremum taken over all intervals $I$ containing $x$.
As is known \cite[p.\,5]{st}, it satisfies the following 
weak-type inequality
\be \lb{5}
   m\{x \in [a,b]: M(f,x) > t\} \le \frac{5}{t}\, \|f\|_1\,.
\ee
\end{definition}

\begin{proposition}
For any $\Delta_n$, we have
\be \lb{M1}
   |\Pn(f,x)| \le c_k M(f,x), \qquad  x \in [a,b]\,.
\ee
\end{proposition}

\proof
Let $x\in [a,b]$, and let the index $i$ be such that $x\in I_i$ 
and $|I_i| \ne 0$. By definition 
of the Dirichlet kernel $K_{\Delta_n}$,
$$
    P_{\Delta_n}(f,x) = \int_{a}^b K_{\Delta_n}(x,y)f(y)\,dy\,,
$$
so using inequality \rf[Kineq] from the previous lemma, we obtain
$$
    |P_{\Delta_n}(f,x)|
\le \sum_{j=1}^n \int_{I_j} |K_{\Delta_n}(x,y)| |f(y)|\,dy
\le   C \sum_{j=1}^n  \frac{\theta^{|i-j|}}{|I_{ij}|}\int_{I_j} |f(y)|\,dy\,.
$$
Since $I_j\subset I_{ij}$ and $x\in I_i\subset I_{ij}$, the definition 
\rf[M0] of the maximal function implies
$\int_{I_j}|f(y)|\,dy \le \int_{I_{ij}}|f(y)|\,dy \le |I_{ij}| M(f,x)$. Hence,
$$
|P_{\Delta_n}(f,x)|\leq C \sum_{j=1}^n \theta^{|i-j|} M(f,x)\,,
$$
and \rf[M1] is proved.
\qed

\medskip
On combining \rf[M1] and \rf[5], we obtain a weak-type
inequality for $P^*$.

\begin{corollary}
For the maximal operator $P^*(f,x) := \sup_n |\Pn(f,x)|$, we have
\be \lb{P1}
   m \{x \in [a,b]: P^*(f,x) > t\} \le \frac{c_k}{t}\, \|f\|_1\,.
\ee
\end{corollary}

The next statement is a straightforward corollary of \rf[dBc];
we give its proof for completeness.

\begin{proposition}
We have  
\be \lb{C1}
     f \in C[a,b] \ra
     P_{\Delta_n}(f,x) \to f(x) \quad \mbox{uniformly.}
\ee
\end{proposition}

\proof  
Since $\Pn$ is a linear projector and $\|\Pn\|_\infty \le c_k$ by \rf[dBc], 
the Lebesgue inequality gives us
$$
   \|f - \Pn f\|_\infty \le (c_k + 1)\,E_{\Delta_n}(f)\,,
$$
where $E_{\Delta_n}(f)$ is the error of the best approximation of $f$ 
by splines from $\SS_k(\Delta_n)$ in the uniform norm. It is known that
$$
    E_{\Delta_n}(f) \le c_k \omega_k(f,|\Delta_n|)\,,
$$
where $\omega_k(f,\delta)$ is the $k$-th modulus of smoothness of $f$.
Since $\omega_k(f,\delta) \to 0$ as $\delta\to 0$, 
we have the uniform convergence
$$
   \|f - \Pn f\|_\infty \to 0 \qquad (|\Delta_n| \to 0)\,,
$$
and that proves \rf[C1].
\qed

\medskip\noindent
{\bf Proof of Theorem \ref{t1}.}
The derivation of the almost everywhere convergence of $\Pn f$ for 
$f\in L_1$ from the weak-type 
inequality \rf[P1] and convergence on the dense subset \rf[C1]
follows a standard scheme which can be found in \cite[pp.\,3-4]{g}. 
We present this argument for completeness.

Let $v\in L^1[a,b]$. We define 
$$
    R(v,x)
 := \limsup_{n\to \infty} \Pn v(x) - \liminf_{n\to \infty} \Pn v(x)
$$
and note that $R(v,x) \le 2 P^*(v,x)$, therefore, by \rf[P1],
\be \lb{Rh}
      m \{x \in [a,b]: R(v,x) > \delta\} \le \frac{2c_k}{\delta}\, \|v\|_1\,.
\ee
Also, for any continuous function $g$ we have $R(g,x) \equiv 0$ by \rf[C1], 
and since $\Pn$ is linear, 
$$
   R(f,x) \le  R(f-g,x) + R(g,x) = R(f-g,x). 
$$
This implies, for a given $f\in L_1$ and any $g \in C$,
\baa
      m\,\{x \in [a,b]: R(f,x) > \delta\} 
\le   m\,\{x \in [a,b]: R(f-g,x) > \delta\} 
\stackrel{\rf[Rh]}{\le} \frac{2c}{\delta}\, \|f-g\|_1\,.
\eaa
Letting $\|f-g\|_1 \to 0$, we obtain, for every $\delta>0$, 
$$
      m\,\{x \in [a,b]: R(f,x) > \delta\} = 0\,,
$$
so $R(f,x)=0$ for almost all $x\in [a,b]$. This means that 
$\Pn f$ converges almost everywhere. It remains to show that this limit 
equals $f$ a.e., but this is obtained by replacing 
$R(f,x)$ by $|\lim\limits_{n\to\infty} P_n f(x)-f(x)|$ in the above argument.
\qed


\section{Convergence of $\Pn(f)$ at the Lebesgue points} \lb{clp}


Here, we show that the estimate \rf[Kineq] for the Dirichlet kernel
implies convergence of $\Pn(f,x)$ at the Lebesgue points of $f$. 
Since by the classical Lebesgue differentiation theorem
the set of all Lebesgue points has the full measure,  
this gives a more precise version of Theorem \ref{t1}. 

We use standard arguments similar to those used in 
\cite[Chapter 1, Theorem 2.4]{dl} for integral operators, or in
\cite[Chapter 5.4]{hw} for wavelet expansions. 
 
Recall that a point $x$ is said to be a Lebesgue point of $f$, if 
$$
     \lim_{I\ni x,\,|I|\to 0}|I|^{-1}\int_I |f(x)-f(y)|\, dy = 0\, ,
$$
where the limit is taken over all intervals $I$ containing the point $x$, as the diameter of $I$ tends to zero.

\begin{theorem}
Let $x$ be a Lebesgue point of the integrable function $f$, and let 
$(\Delta_n)$ be a sequence of partitions of $[a,b]$ with $|\Delta_n|\to 0$. 
Then,
\[
   \lim_{n \to \infty} P_{\Delta_n}(f,x) = f(x).
\]
\end{theorem}

\proof
Let $x$ be a Lebesgue point of $f$. Since the spline space $\SS_k(\Delta_n)$
contains constant functions, we have
$\int_a^b K_{\Delta_n}(x,y)\,dy = 1$ for any $x \in [a,b]$, so we need
to prove that
\be \lb{Kto0}
     \int_a^b K_{\Delta_n}(x,y)[f(x) - f(y)]\,dy \to 0 \quad (n \to  \infty).
\ee
For $r > 0$, set $B_r (x) := [x-r, x+r] \cap [a,b]$. 
Now, given $\e > 0$, let $\delta$ be such that 
\be \lb{x}
   |I|^{-1} \int_I |f(x)-f(y)|\,dy < \e \qquad     
\ee
for all intervals $I$ with $I \subset B_{2\delta}(x)$ and $I \ni x\,.$
Further, with $\theta \in (0,1)$ from inequality \rf[Kineq],
take $m$ and $N = N(m)$ such that 
$$
     \theta^m < \e\delta\,,\qquad
      (m+2)|\Delta_n| < \delta \quad \forall n \ge N\,,
$$
and consider any such $\Delta_n$. 

1) Let $|x-y| > \delta$, and let $x \in I_i$ and $y \in I_j$.
Then $|i-j| > m$ and $|I_{ij}| > \delta$, hence, by inequality \rf[Kineq] 
for the Dirichlet kernel $K_{\Delta_n}$,
\[
   |K_{\Delta_n}(x,y)|\leq C\theta^m\delta^{-1}\leq C\e\,.
\]
As a consequence,
\be \lb{greaterdelta}
     \int_{|x-y|>\delta} |K_{\Delta_n}(x,y)||f(x)-f(y)|\,dy
\leq C\e \int_a^b |f(x)-f(y)|\, dy \le 2C \e \|f\|_1\,.
\ee

2)  Let $|x-y| \le \delta$, i.e.\ $y \in B_\delta(x)$, and let $x \in I_i$.
Note that if 
$I_j \cap B_\delta(x) \ne \emptyset$, then $I_j \subset B_{2\delta}(x)$,
hence $I_{ij} \subset B_{2\delta}(x)$ as well, and again, 
by inequality \rf[Kineq],
\baa
      \int_{B_\delta(x)}|K_{\Delta_n}(x,y)||f(x)-f(y)|\, dy 
&\le& \sum_{j:I_j\cap B_\delta(x)\neq\emptyset} 
        \int_{I_j} |K_{\Delta_n}(x,y)||f(x)-f(y)|\, dy \\
&\le& C \sum_{j:I_{ij}\subset B_{2\delta}(x)}\theta^{|i-j|} \Big(|I_{ij}|^{-1} 
        \int_{I_{ij}} |f(x)-f(y)|\,dy \Big).
\eaa
By \rf[x], since $x \in I_{ij}\subset B_{2\delta}(x)$, 
the terms in the parentheses are all bounded by $\e$, 
therefore
\be \lb{smallerdelta}
     \int_{|x-y|<\delta}|K_{\Delta_n}(x,y)||f(x)-f(y)|\, dy  
\le C\e \sum_{j} \theta^{|i-j|} \le C_1 \e\,.
\ee
Combining estimates \rf[greaterdelta] and \rf[smallerdelta] 
for the integration over $|x-y|>\delta$ and $|x-y|<\delta$, respectively, 
we obtain \rf[Kto0], i.e. convergence of $P_{\Delta_n}(f,x)$ 
to $f(x)$ at Lebesgue points of $f$, provided $|\Delta_n|\to 0$.
\qed


\section{Proof of Theorem \ref{p1}} \lb{pp1}


We will prove \rf[a0] for $i \le j$. This proves also the case
$i \ge j$,  since $h_{ij} = h_{ji}$ and $a_{ij} = a_{ji}$. 

So, for the entries $\{a_{ij}\}$ of the inverse of the matrix 
$\{\langle N_i,N_j\rangle\}$, we want to show that 
\be \lb{aa0}
    |a_{ij}| \le  K\, \gamma^{|i-j|} h_{ij}^{-1}\,,
\ee
where $h_{ij}$ is the length of a largest subinterval of $[t_i,t_{j+k}]$.
The proof is based on the following lemma.

\begin{lemma} \lb{le1}
For any $\Delta_n$, let $\{a_{ij}\}$ be the inverse of the B-spline Gram 
matrix $\{\langle N_i,N_j\rangle\}$. Then
\ba
    |a_{is}| 
&\le& K_1\, \gamma^{|i-s|} (\max\{\kappa_i, \kappa_s\})^{-1}\,, \lb{a1} \\
    |a_{ij}| 
&\le& K_2\, \gamma^{|\l-j|} 
        \sum_{\mu=\l-(k-1)}^{\l+k-2} |a_{i\mu}|\,, \qquad
        i + k \le \l < j\,,\lb{a2} \\
    |a_{i\mu}| 
&\le& K_3\, \max_{\mu-(k-1) \le s \le \mu-1} |a_{is}|\,, \qquad 
     i < \mu\,, \lb{a3}
\ea
where $K_i > 0$ and $\gamma \in (0,1)$ are some constants that 
depend only on $k$.
\end{lemma}

\begin{remark} \rm
All three estimates are known in a sense.
Inequalities \rf[a1] and \rf[a2] follow from Demko's theorem \cite{d}
on inverses of band matrices and the fact \cite{s2} 
that the inverse of the Gram matrix $G = \{(M_i,N_j)\}_{i,j=1}^n$
satisfies $\|G^{-1}\|_\infty < c_k$.
Actually, \rf[a1] was explicitly given by Ciesielski \cite{c}, 
while \rf[a2] is a part of Demko's proof. Inequality \rf[a3]
appeared in Shadrin's manuscript \cite{s1}, and it does not
use the uniform boundedness of $\|G^{-1}\|_\infty$. As those estimates
are scattered in the aforementioned papers, we extract the relevant parts from them and present the proofs of \rf[a1]--\rf[a3] in Sect.\,\ref{pl1}.
\end{remark}

\medskip\noindent
{\bf Proof of Theorem \ref{p1}.} 
Let $I_\l$ be a largest subinterval of $[t_i,t_{j+k}]$, i.e.,
$$
   h_{ij} = \max\, \{h_s\}_{s=i}^{j+k-1} = h_\l\,.
$$

1) If $I_\l$ belongs to the support of $N_i$ or that of $N_j$, then 
$$
   \max(\kappa_i,\kappa_j) \ge h_\l = h_{ij}\,,
$$
and, by \rf[a1],
$$
    |a_{ij}| 
\le K_1\, \gamma^{|i-j|} (\max\{\kappa_i, \kappa_j\})^{-1}\
\le  K_1\, \gamma^{|i-j|} h_{ij}^{-1}\,,
$$
so \rf[aa0] is true. 

\medskip
2) Now, assume that $I_\l$ does not belong to the supports of 
either $N_i$ or $N_j$, i.e.,
$$
    i + k \le \l < j\,.
$$
Consider the B-splines $(N_s)_{s=\l+1-k}^\l$ whose support 
$[t_s,t_{s+k}]$ contains $I_\l = [t_\l,t_{\l+1}]$. Then 
$$
    \kappa_s \ge h_\l = h_{ij}\,, \qquad \l-(k-1) \le s \le \l\,.
$$
Using estimate \rf[a1], we obtain for such $s$
$$
    |a_{is}| 
\le K_1\, \gamma^{|i-s|} \kappa_s^{-1}
\le K_1\, \gamma^{|i-s|} h_{ij}^{-1}
\le K_1\, \gamma^{-k} \gamma^{|i-\l|} h_{ij}^{-1}, 
$$
i.e.,
\be \lb{a4}
    \max_{\l-(k-1) \le s \le \l} |a_{is}| 
\le C_1\, \gamma^{|i-\l|} h_{ij}^{-1}\,. 
\ee

3) From \rf[a2], we have
\be \lb{a5}
    |a_{ij}| 
\le  2(k-1) K_2 \gamma^{|\l-j|} 
     \max_{\l-(k-1) \le \mu \le \l+k-2} |a_{i\mu}|\,.
\ee
Note that \rf[a3] bounds $|a_{i\mu}|$ in terms of of the absolute values 
of the $k-1$ coefficients that precede it, hence by induction and with the
understanding that $K_3 > 1$, 
$$
    |a_{i,\l+r}| 
\le K_3^r \max_{\l-(k-1) \le s \le \l} |a_{is}|, \qquad r=1,2,\ldots,
$$
therefore
\be \lb{a6}
     \max_{\l-(k-1) \le \mu \le \l+k-2} |a_{i\mu}|
\le K_3^{k-2} \max_{\l-(k-1) \le s \le \l} |a_{is}|.
\ee
Combining \rf[a5], \rf[a6] and \rf[a4], gives
$$
    |a_{ij}| 
\le 2(k-1)K_2 \gamma^{|\l-j|} K_3^{k-2} C_1 \gamma^{|i-\l|} h_{ij}^{-1}
 =  K\, \gamma^{|i-j|} h_{ij}^{-1}\,,
$$
and that proves \rf[aa0], hence \rf[a0]. 
\qed


\section{Proof of Lemma \ref{le1}} \lb{pl1}


Here, we prove the three parts of Lemma \ref{le1} as Lemmas \ref{le11},
\ref{le12} and \ref{le13}, respectively. 
The proof is based on certain properties of the Gram matrix 
$G := \{\la M_i,N_j\r\}_{i,j=1}^n$ and its inverse 
$G^{-1} =: \{b_{ij}\}_{i,j=1}^n$. Here, $(M_i)$ is the sequence
of the $L_1$-normalized B-splines on $\Delta_n$, 
$$
     M_i := \frac{k}{\kappa_i}\,N_i\,, \qquad 
     \int_{t_i}^{t_{i+k}} M_i(t) \,dt = 1.
$$

First, we note that $G$ is a banded matrix with max-norm one, i.e.,
\be \lb{G1}
     \la M_i, N_j\r = 0 \quad\mbox{for}\quad |i-j| > k-1, \qquad
    \|G\|_\infty = 1,
\ee
where the latter equality holds due to the fact that
$\sum_j |\la M_i, N_j\r| = \la M_i, \sum_j N_j\r = \la M_i,1\r = 1$.
A less obvious property is the boundedness of $\|G^{-1}\|_\infty$. 

\begin{theorem}[Shadrin \cite{s2}]
For any $\Delta_n$, with $G := \{\la M_i,N_j\r\}_{i,j=1}^n$, we have 
\be \lb{G2}
     \|G^{-1}\|_\infty \le c_k\,,
\ee
where $c_k$ is a constant that depends only on $k$.
\end{theorem}
We recall that \rf[G2] is equivalent to \rf[dBc], i.e., the $\l_\infty$-norm 
boundedness of the inverse $G^{-1}$ of the Gramian is equivalent to 
the $L_\infty$-norm boundedness of the orthogonal spline 
projector $\Pn$, namely, with some constant $d_k$ (e.g., 
the same as in \rf[d_k]), we have
$$
    \frac{1}{d_k^2}\, \|G^{-1}\|_\infty 
\le \|\Pn\|_\infty \le \|G^{-1}\|_\infty\,.
$$

Next, we apply the following theorem to $G$. 

\begin{theorem}[Demko \cite{d}] \lb{Demko}
Let $A=(\alpha_{ij})$ be an $r$-banded matrix, i.e., 
$\alpha_{ij} = 0$ for $|i-j| > r$,
and let $\|A\|_p \le c'$ and  $\|A^{-1}\|_p \le c''$ for some 
$p \in [1,\infty]$.  Then the elements of the inverse 
$A^{-1} =: (\alpha_{ij}^{(-1)})$ decay exponentially away from the diagonal, 
precisely
$$
    |\alpha_{ij}^{(-1)}| \le K \gamma^{|i-j|}\,, 
$$
where $K > 0$ and $\gamma \in (0,1)$ are constants that depend only on 
$c'$, $c''$ and $r$.
\end{theorem}
We will need two corollaries of this result.
\begin{corollary} \lb{coK}
For any $\Delta_n$, with $G = \{\la M_i,N_j\r\}_{i,j=1}^n$, 
and $G^{-1} =: \{b_{ij}\}_{i,j=1}^n$,
we have
\be \lb{b}
    |b_{ij}| \le K_0 \gamma^{|i-j|}\,,
\ee
where $K_0 > 0$ and $\gamma \in (0,1)$ are constants that depend only on $k$.
\end{corollary}

\proof
Indeed, by \rf[G1]-\rf[G2], we may apply Demko's theorem to the Gram 
matrix $G$, with $c'=1$, $c'' = c_k$, $r =k-1$, and $p=\infty$, 
and that gives the statement.
\qed

\begin{corollary} 
For any $\Delta_n$, with $G = \{\la M_i,N_j\r\}_{i,j=1}^n$, 
\be \lb{c1}
     \|G\|_1 < c_1, \qquad\|G^{-1}\|_1 < c_2\,,
\ee
where $c_1,c_2$ depend only on $k$.
\end{corollary}

\proof
It follows from \rf[b] that $\|G^{-1}\|_1 = \max_j \sum_i |b_{ij}|$ 
is bounded, whereas $\|G\|_1$ is bounded since $G$ is a $(k-1)$-banded 
matrix with nonnegative entries $\la M_i,N_j\r \le 1$.
\qed




\medskip
Now we turn to the proof of Lemma \ref{le1} starting with inequality \rf[a1]. 

\begin{lemma}[{Property \rf[a1]}] \lb{le11}
Let $\{a_{ij}\}$ be the inverse of the B-spline Gram matrix 
$\{\langle N_i,N_j\rangle\}$.
Then
\be \lb{aa1}
    |a_{is}| 
\le K_1\, \gamma^{|i-s|} (\max\{\kappa_i,\kappa_s\})^{-1}\,.
\ee
\end{lemma}

\proof
As we mentioned earlier, this estimate was proved by Ciesielski 
\cite[Property 6]{c}.  Here are the arguments. 
The elements of the two inverses $\{a_{ij}\} = G_0^{-1}$ and 
$\{b_{ij}\} = G^{-1}$ are connected by the formula
\be \lb{ab}
    a_{ij} = b_{ij}(k/\kappa_j) = b_{ji}(k/\kappa_i)\,.
\ee
Indeed, the identity $N_i = \kappa_i M_i /k $ implies that
the matrix $G_0 := \{\la N_i,N_j\r\}$ is related to $G = \{\la M_i,N_j\r\}$ in the form
$$
   G_0 = D G, \qquad \text{where } D = {\rm diag}\, [\kappa_1/k,\ldots,\kappa_n/k]\,.
$$
Hence, $G_0^{-1} =  G^{-1} D^{-1}$, and the first equality in \rf[ab] follows.
The second equality is a consequence of the symmetry of $G_0$, 
as then $G_0^{-1}$ is symmetric too, i.e., $a_{ij} = a_{ji}$.
Then, in \rf[ab] we may use the estimate $|b_{ij}| \le K_0 \gamma^{|i-j|}$ 
from \rf[b], and \rf[aa1] follows.
\qed




\begin{lemma}[{Property \rf[a2]}] \lb{le12}
Let $\{a_{ij}\}$ be the inverse of the B-spline Gram matrix 
$\{\langle N_i,N_j\rangle\}$.
Then
\be \lb{aa2}
    |a_{ij}| 
\le K_2\, \gamma^{|\l-j|} \sum_{\mu=\l-(k-1)}^{\l+k-2} |a_{i\mu}|\,,
     \qquad i+k \le \l < j\,.
\ee
\end{lemma}

\proof
1) Since $a_{ij} =  b_{ji}(k/\kappa_i)$ by \rf[ab], it is sufficient to 
establish the same inequality for the elements $b_{ji}$
of the matrix $G^{-1} = (b_{ij})$:
\be \lb{bb1}
    |b_{ji}| 
\le K_2\, \gamma^{|\l-j|} \sum_{\mu=\l-(k-1)}^{\l+k-2} |b_{\mu i}|\,.
\ee
We fix $i$ with $1 \le i \le n$, and to simplify notations
we write $b_j := b_{ji}$, omitting $i$ in the subscripts. So,
the vector $b = (b_1,\ldots,b_n)^T$ is the $i$-th 
column of $G^{-1}$, hence
\be \lb{e}
     G b = e_i\,.
\ee

2) The following arguments just repeat those
in the proof of Theorem \ref{Demko} used by Demko \cite{d} 
and extended by de\,Boor \cite{b}.

For $m > i$, set 
$$
    b^{(m)} = (0,0,\ldots,0,b_m, b_{m+1}, \ldots b_n)^T\,.
$$
With $r := k-1$, the Gram matrix $G = \{(M_i,N_j)\}_{i,j=1}^n$ 
is $r$-banded, and that together with \rf[e] implies 
$$
    \supp G\,b^{(m)} \subset [m-r, m+(r-1)]\,.
$$
It follows that $G\,b^{(m)}$ and $G\,b^{(m+2r)}$ have disjoint support,
therefore
$$
     \|G\,b^{(m)}\|_1 + \|G\,b^{(m+2r)}\|_1 
=  \|G\,b^{(m)} - G\,b^{(m+2r)}\|_1\,.
$$
This yields
\baa
      \|G^{-1}\|_1^{-1}\, (\,\|b^{(m)}\|_1 + \|b^{(m+2r)}\|_1)
&\le&  \|G\,b^{(m)}\|_1 + \|G\,b^{(m+2r)}\|_1 \\
& = &  \|G\,b^{(m)} - G\,b^{(m+2r)}\|_1 \\
&\le&  \|G\|_1\, \|b^{(m)} - b^{(m+2r)}\|_1 \\
& = &  \|G\|_1\, (\,\|b^{(m)}\|_1 - \|b^{(m+2r)}\|_1)\,,
\eaa
i.e.,
\be \lb{c3}
    \|b^{(m)}\|_1 + \|b^{(m+2r)}\|_1
\le c_3\, (\,\|b^{(m)}\|_1 - \|b^{(m+2r)}\|_1)\,,
\ee
where $c_3 = \|G\|_1 \|G^{-1}\|_1 \ge 1$. 
This gives
$$
   \|b^{(m+2r)}\|_1 \le \gamma_0\,  \|b^{(m)}\|_1, \qquad
   \gamma_0 = \frac{c_3-1}{c_3 + 1} < 1\,,
$$
where $\gamma_0$ depends only on $k$ since so does $c_3 = c_1c_2$ by \rf[c1]. 

It follows that, for any $j,m$ such that $i < m  \le j \le n$, 
we have
\baa
    |b_j| 
&\le& \|b^{(j)}\|_1 
\le \gamma_0^{\lfloor \frac{j-m}{2r} \rfloor}\,\|b^{(m)}\|_1
\le \gamma_0^{-1}\gamma_0^{|j-m|/2r}\,\|b^{(m)}\|_1 \\
&=:& c_4\, \gamma^{|j-m|}\,\|b^{(m)}\|_1\,.
\eaa
Applying \rf[c3] to the last line, we obtain
$$
    |b_j| 
\le c_4\, \gamma^{|j-m|}\, c_3\, (\|b^{(m)}\|_1 - \|b^{(m+2r)}\|_1) 
= c_5\, \gamma^{|j-m|} \sum_{\mu=m}^{m+2r-1} |b_\mu| \,.
$$
Taking $m = \l - r = \l - (k-1)$, we bring this inequality 
to the form \rf[bb1] needed:
$$
    |b_j| 
= c_5\,\gamma^{k-1} \gamma^{|j-\l|} \sum_{\mu=\l-(k-1)}^{\l+k-2} |b_\mu| \,.
$$




\begin{lemma}[{Property \rf[a3]}] \lb{le13}
Let $\{a_{ij}\}$ be the inverse of the B-spline Gram matrix 
$\{\langle N_i,N_j\rangle\}$.
Then
\be \lb{aa3}
    |a_{im}| 
\le K_3\, \max_{m-(k-1) \le s \le m-1} |a_{is}|\,, \qquad m > i,
\ee
i.e., the absolute value of a coefficient following $a_{ii}$ can be bounded
in terms of the absolute values of the $k-1$ coefficients directly preceding 
that coefficient.
\end{lemma}

\proof
This estimate appeared in \cite[proof of Lemma 7.1]{s1}.
To adjust that proof to our notations, we note that the basis $\{N_i^*\}$
dual to the B-spline basis $\{N_i\}$ is given by the formula 
$$
   N_i^* = \sum_{j=1}^n a_{ij} N_j\,.
$$
Indeed, from the definition of $a_{ij}$,  
we have $\la N_i^*,N_m \r = \sum_{j=1}^n a_{ij} \la N_j, N_m\r = \delta_{im}$.

1) We fix $i$, write $a_j := a_{ij}$ omitting the index $i$, 
and for $m > i$, set 
\be \lb{psi0}
    \psi_{m-(k-1)} := \sum_{j=m-(k-1)}^n a_j N_j\,, \qquad
    \psi_{m} := \sum_{j=m}^n a_j N_j\,.
\ee
Then, since $\supp N_j = [t_j,t_{j+k}]$, it follows that 
$$
    \psi_{m-(k-1)}(x) = N_i^*(x), \qquad x \in [t_m,b]\,,
$$
Therefore, $\psi_{m-(k-1)}$ is orthogonal to $\spa\{N_j\}_{j = m}^n$,
in particular to $\psi_m$. This gives
\be \lb{psi1}
   \|\psi_{m-(k-1)}\|_{L_2[t_m,b]}^2 + \|\psi_m\|^2_{L_2[t_m,b]}
= \|\psi_{m-(k-1)} - \psi_m\|_{L_2[t_m,b]}^2\,.
\ee

2) Further, we have
$$
   E_m = [t_m,t_{m+k}] \subset [t_m,b]\,,
$$
whereas the equality 
$\psi_{m-(k-1)} - \psi_m = \sum_{j=m-(k-1)}^{m-1} a_j N_j$
implies
$$
   \supp (\psi_{m-(k-1)} - \psi_m) \cap [t_m,b]
 = [t_m,t_{m+k-1}] \subset E_m\,.
$$
Therefore, from \rf[psi1], we conclude
\be \lb{psi}
    \|\psi_{m-(k-1)}\|_{L_2(E_m)}^2 + \|\psi_m\|^2_{L_2(E_m)}
\le \|\psi_{m-(k-1)} - \psi_m\|_{L_2(E_m)}^2\,.
\ee

3) Now recall that, by a theorem of de Boor (see \cite{b1} 
or \cite[Chapter 5, Lemma 4.1]{dl}), there is a constant
$d_k$ that depends only on $k$ such that 
\be \lb{d_k}
    d_k^{-2} |c_m|^2
\le  |E_m|^{-1} \|\sum_{j=1}^n c_j N_j\|_{L_2(E_m)}^2 \qquad \forall c_j \in \R\,.
\ee
(This gives the upper bound $d_k$ for the B-spline basis condition number.)
So, applying this estimate to the left-hand side of \rf[psi], where we use 
\rf[psi0], we derive
\baa
    2\, d_k^{-2} |a_m|^2
&\le& |E_m|^{-1} \left( \|\psi_{m-(k-1)}\|_{L_2(E_m)}^2 
         + \|\psi_m\|^2_{L_2(E_m)} \right) \\
&\stackrel{\rf[psi]}{\le}& |E_m|^{-1}\,\|\psi_{m-(k-1)} - \psi_m\|_{L_2(E_m)}^2 \\
&\le& \|\psi_{m-(k-1)} - \psi_m\|_{L_\infty(E_m)}^2 \\
&\stackrel{\rf[psi0]}{=}& \| \sum_{j=m-(k-1)}^{m-1} a_j N_j \|_{L_\infty(E_m)}^2 \\
&\le& \max_{m-(k-1) \le s \le m-1} |a_s|^2\,,
\eaa
i.e.,
$$
    |a_m| \le K_3\,\max_{m-(k-1) \le s \le m-1} |a_s|^2\,,
$$
and that proves \rf[aa3].
\qed

\medskip
{\bf Acknowledgements.} 
The authors are grateful to Carl de Boor for his help in arranging
this collaboration and for many valuable remarks.  
It is also a pleasure to thank Zbigniew Ciesielski, 
Anna Kamont and Paul M\"uller for their comments on a draft of this paper. 
Finally, we thank the anonymous referee for various remarks, in particular 
for pointing out convergence at Lebesgue points.

The first author is supported by the Austrian Science Fund, FWF project P\,23987-N18, and a part of this research was performed while he was visiting the Institute of Mathematics of the Polish Academy of Sciences 
in Sopot. He thanks the Institute  for its hospitality and 
excellent working conditions. These stays were supported by 
MNiSW grant N N201 607840.


\medskip
\textsc{Markus Passenbrunner, Institute of Analysis, 
Johannes Kepler University Linz, Altenberger Strasse 69, 4040 Linz, Austria} \\
\emph{E-mail address: } \texttt{markus.passenbrunner@jku.at}

\medskip
\textsc{Alexei Shadrin, DAMTP, University of Cambridge, Wilberforce Road, 
Cambridge CB3 0WA, UK} \\
\emph{E-mail address: } \texttt{a.shadrin@damtp.cam.ac.uk}

\end{document}